\documentclass[11pt,reqno]{amsart}
\usepackage[T1]{fontenc}
\usepackage[utf8]{inputenc}
\usepackage{lmodern}
\usepackage[a4paper,margin=28mm,headheight=14pt]{geometry}
\usepackage{amsmath,amssymb,mathtools}
\usepackage{microtype}
\usepackage{enumitem}
\usepackage{needspace}
\usepackage[hidelinks]{hyperref}
\hypersetup{
  pdftitle={Stability and exact extremal colourings for the Alon--Frankl--Lovasz bound},
  pdfsubject={Stability and eventual exact classification of monochromatic matching extremizers},
  pdfkeywords={hypergraph, monochromatic matching, stability, extremal colouring, weighted defect}
}
\newtheorem{theorem}{Theorem}[section]
\newtheorem{lemma}[theorem]{Lemma}

\theoremstyle{definition}
\newtheorem{definition}[theorem]{Definition}
\theoremstyle{remark}

\numberwithin{equation}{section}
\DeclareMathOperator{\supp}{supp}
\DeclareMathOperator{\conv}{conv}
\DeclareMathOperator{\dist}{dist}

\newcommand{\Prob}{\mathbb P}
\newcommand{\one}{\mathbf 1}
\newcommand{\ee}{\mathbf e}
\newcommand{\cC}{\mathcal C}
\newcommand{\cE}{\mathcal E}

\newcommand{\eps}{\varepsilon}
\newcommand{\norminf}[1]{\left\lVert #1\right\rVert_\infty}
\newcommand{\floor}[1]{\left\lfloor #1\right\rfloor}
\newcommand{\ceil}[1]{\left\lceil #1\right\rceil}

\title[Stability and extremal colourings for the AFL bound]{Stability and exact extremal colourings for the\newline
Alon--Frankl--Lov\'asz bound}
 \author{Tianhong Kan, Hongliang Lu and Feihong Yuan}
 \address{School of Mathematics and Statistics, Xi'an Jiaotong University}
 \email{kantianhong@outlook.com;luhongliang@mail.xjtu.edu.cn; fhyuan1@gmail.com}

\begin{document}
\begin{abstract}
	The Alon--Frankl--Lov\'asz theorem states that every $q$-colouring
	of the edges of the complete $r$-uniform hypergraph on $n$
	vertices contains a monochromatic matching of size
	$\lfloor(n+q-1)/(r+q-1)\rfloor$. For fixed uniformity $r$ and number
	of colours $q$, we prove a stability theorem: every colouring
	whose largest monochromatic matching has size at most
	$n/(r+q-1)+o(n)$ differs from a threshold extremal colouring
	on only $o(n^r)$ edges. We also give a necessary and sufficient
	structural characterization of all extremal colourings for
	sufficiently large $n$. The proofs combine rigidity of fractional covering weights
	with an extension argument for prescribed matchings. These results answer, for sufficiently large orders, the
	stability and extremal-colouring questions posed by
	Gishboliner, Glock, Michaeli and Sgueglia.
\end{abstract}
\maketitle

\section{Introduction}\label{sec:intro}

Packing vertex-disjoint monochromatic subgraphs is a natural
theme in Ramsey theory. One instance is the search
for a monochromatic matching, a collection of pairwise disjoint
edges of the same colour. The size of such a matching depends
strongly on how the edges of that colour are distributed:
even a large family of edges may contain only small matchings
if every edge meets a fixed small set of vertices. In a colouring
of a complete graph or hypergraph, the colour classes partition
the entire edge set, so restrictions on their matching numbers
must be satisfied simultaneously. This interaction between
disjointness within each colour and the requirement to colour
every edge makes monochromatic matchings a natural meeting
point of Ramsey theory and extremal set theory.

For a positive integer $a$, write $[a]=\{1,\ldots,a\}$, and let
$\binom Vr$ denote the family of $r$-element subsets of $V$.
An $r$-graph is an $r$-uniform hypergraph, and $K_n^{(r)}$ is the
complete $r$-graph on $n$ vertices. A matching is a set of pairwise vertex-disjoint edges; its size
is its number of edges. Write $\nu(H)$ for the maximum size
of a matching in $H$. For a colouring
$\chi:\binom Vr\to[q]$, let $H_i=H_i(\chi)$ be the spanning
$r$-graph with edge set $\chi^{-1}(i)$.
Throughout, $r,q\ge2$ are fixed and $n=|V|\ge r$. Put
\begin{equation}\label{eq:parameters}
 \begin{gathered}
 L=r+q-1,\qquad k=\floor{\frac{n+q-1}{L}},\qquad
 m=k+1,\qquad D=Lm-q-n.
 \end{gathered}
\end{equation}
The definition of $k$ gives
$Lk-q+1\le n\le Lk+r-1,\ 0\le D\le L-1$.
The Alon--Frankl--Lov\'asz theorem~\cite{AFL} asserts that
\begin{equation}\label{eq:AFL}
 \max_{i\in[q]}\nu(H_i)\ge k.
\end{equation}
This bound is sharp for every $n\ge r$. We call a colouring \emph{extremal} if
equality holds, equivalently if $\nu(H_i)\le k$ for every $i$.

The graph case $r=2$ is the equal-size case of the matching
Ramsey theorem of Cockayne and Lorimer~\cite{CL}. More generally,
the bound has a natural Kneser-hypergraph interpretation: if
$\operatorname{KG}^{(m)}(n,r)$ has vertex set $\binom{[n]}r$ and
hyperedges consisting of $m$ pairwise disjoint $r$-sets, then a
$q$-colouring of $K_n^{(r)}$ with no monochromatic matching of size
$m$ is precisely a proper $q$-colouring of
$\operatorname{KG}^{(m)}(n,r)$. The case $m=2$ goes back to the
proofs of Kneser's conjecture by Lov\'asz~\cite{Lovasz78} and
B\'ar\'any~\cite{Barany78}. The general result of
Alon, Frankl and Lov\'asz~\cite{AFL} belongs to this topological
tradition; see Matou\v{s}ek~\cite{Matousek} for background.

Gishboliner, Glock, Michaeli and Sgueglia~\cite[Section~4]{GGMS}
asked for a stability theorem and a characterization of the
extremal colourings. We prove the stability statement and obtain
an exact characterization for all sufficiently large $n$, in
every residue class modulo $L$.

\subsection{Threshold colourings}
The extremal constructions are described by positive integer
thresholds whose sum is $L$. Set
\begin{equation*}\label{eq:thresholds}
 \cC_{r,q}=\left\{c\in\mathbb Z_{>0}^q:\sum_i c_i=L\right\}.
\end{equation*}
Every such vector satisfies
\begin{equation}\label{eq:threshold-range}
 1\le c_i\le r,\qquad \sum_i(c_i-1)=r-1.
\end{equation}
The upper bound follows from $c_j\ge1$ for the other $q-1$
coordinates.

\begin{definition}\label{def:threshold}
Let $\cE_{r,q}(V)$ be the family of colourings $\psi$ for which
there are $c\in\cC_{r,q}$ and a partition into possibly empty parts
$V=V_1\sqcup\cdots\sqcup V_q$ such that
\begin{equation*}\label{eq:threshold-rule}
 |V_i|\le c_i m-1,\qquad
 \psi(e)=i\ \Longrightarrow\ |e\cap V_i|\ge c_i.
\end{equation*}
A colour $i$ with $|e\cap V_i|\ge c_i$ is \emph{eligible} for $e$.
Each edge may be assigned any of its eligible colours; no
compatibility condition is imposed on different edges.
\end{definition}

For each $c\in\cC_{r,q}$, such colourings exist and are extremal. Indeed, the capacities
$u_i=c_i m-1$ are nonnegative and have sum
$Lm-q=n+D\ge n$, so a partition respecting them exists.
Every $r$-edge has an eligible colour, since otherwise its size
would be at most $\sum_i(c_i-1)=r-1$. Finally,
\[
 \nu(H_i(\psi))\le\floor{\frac{|V_i|}{c_i}}\le m-1=k,
\]
and~\eqref{eq:AFL} gives equality for the maximum over colours.
For any partition in Definition~\ref{def:threshold},
\begin{equation}\label{eq:part-sizes}
 0\le u_i-|V_i|\le D,\qquad
 |V_i|=\frac{c_i n}{L}+O_{r,q}(1).
\end{equation}
The first assertion follows because the nonnegative capacity
deficits sum to $D$; the second follows from
$m=n/L+O_{r,q}(1)$.

\subsection{The two results}
For two colourings $\chi,\psi$ on the same vertex set $V$, define
\[
 \dist(\chi,\psi)=|\{e\in\tbinom Vr:\chi(e)\ne\psi(e)\}|,
 \qquad
 \dist(\chi,\cE_{r,q}(V))
 =\min_{\psi\in\cE_{r,q}(V)}\dist(\chi,\psi).
\]

\begin{theorem}[Edge-edit stability]\label{thm:stability}
For every $\eps>0$, there are $\mu=\mu(r,q,\eps)>0$ and
$n_0=n_0(r,q,\eps)$ such that every
colouring $\chi:\binom Vr\to[q]$ with $n\ge n_0$ and
\[
 \nu(H_i)\le(1+\mu)\frac nL\qquad(i\in[q])
\]
satisfies $\dist(\chi,\cE_{r,q}(V))\le\eps n^r$.
\end{theorem}

Equivalently, for a sequence with $n\to\infty$,
\[
 \max_i\nu(H_i)\le n/L+o(n)
 \quad\Longrightarrow\quad
 \dist(\chi,\cE_{r,q}(V))=o(n^r).
\]
The stability theorem permits $o(n^r)$ exceptional edges. Our
second result describes exactly which exceptions are possible
when $\max_i\nu(H_i)=k$. The quantity $D$ is the total unused
capacity of a threshold partition, and the integers $s_i$ below
distribute this unused capacity among the colours.

\Needspace{12\baselineskip}
\begin{theorem}[Exact normal form]\label{thm:exact}
There is $n_{\mathrm{ex}}=n_{\mathrm{ex}}(r,q)$ such that, for
every $n\ge n_{\mathrm{ex}}$, a colouring $\chi$ is extremal
if and only if there are $c\in\cC_{r,q}$, nonnegative integers
$s_1,\ldots,s_q$, and a partition $V=W_1\sqcup\cdots\sqcup W_q$
with
\begin{equation}\label{eq:normal-sizes}
 |W_i|=c_i m-1-s_i,\qquad \sum_i s_i=D,
\end{equation}
such that, for each $i\in[q]$, the family of \emph{violating} edges
\[
 \mathcal F_i=\{e\in\tbinom Vr:\chi(e)=i,\ |e\cap W_i|<c_i\}
\]
satisfies
\begin{equation}\label{eq:budget}
 \sum_{e\in M}\bigl(c_i-|e\cap W_i|\bigr)\le s_i
 \quad\text{for every matching }M\subseteq\mathcal F_i.
\end{equation}
For colour $i$, it suffices to check matchings in $\mathcal F_i$
with at most $s_i+1\le L$ edges.
The sufficiency holds for every $n\ge r$ for which the stated
partition exists.
\end{theorem}

The positive integer $\delta_i(e)=c_i-|e\cap W_i|$ is the
\emph{deficit} of a violating colour-$i$ edge. Thus the theorem
places a budget on the total deficit of disjoint violating edges.
The number of edges in each test is bounded independently of $n$;
the total number of violating edges need not be bounded.
Neither the threshold vector nor the partition is asserted to
be unique.


\subsection{Related work and proof strategy}
The extremal-colouring problem is related to, but distinct from,
the Erd\H{o}s matching problem~\cite{Erdos65}: the latter asks for
the maximum size of a single $r$-uniform family with prescribed
matching number, whereas here all colour classes partition
$\binom Vr$ and must simultaneously have small matching number.
Important bounds for the single-family problem were obtained by
Frankl~\cite{Frankl13} and Frankl and Kupavskii~\cite{FK22}.
Their two basic constructions, a family meeting a fixed small
set and a complete $r$-graph on a restricted vertex set, correspond
to the endpoint conditions $c_i=1$ and $c_i=r$ in our threshold
rule. More general families defined by a lower bound on the
intersection with a fixed set also feature in the complete
intersection theorem of Ahlswede and Khachatrian~\cite{AK} and
the union-restricted families studied by Frankl and
Kupavskii~\cite{FK21}. These analogies explain the relevance of
intersection thresholds, but do not by themselves classify
extremal colourings.

Stability results for bounded matching number provide further
context. In particular, Frankl, Lu, Ma and Wu~\cite{FLMW}
proved a stability theorem for $4$-uniform hypergraphs as part
of their work on the Erd\H{o}s matching conjecture. Their
single-family extremal problem differs from the simultaneous
multicolour problem considered here.

The proof has two stages. In the stability stage, weak regularity
and fractional matching duality associate a vector of covering
weights with each cluster. Near equality in the AFL bound forces
almost all these vectors to lie near the axis points
$\mathbf e_i/c_i$, with proportions $c_i/L$. These clusters give
the required vertex partition; a small adjustment enforces its
exact capacity bounds.

The exact stage begins with this approximate partition. Removing
a small set of exceptional vertices leaves cores in which the
number of edges of an incorrect colour through any vertex is
$o(n^{r-1})$. Each exceptional vertex is assigned a colour that
occurs on a positive proportion of its baseline links to the
cores. A switching argument then completes any prescribed
bounded matching whenever the relevant capacity inequality
fails. Applying this argument first to the empty matching gives
the exact part capacities, and applying it to matchings of
violating edges gives the deficit budgets.

The threshold constructions already appear in the discussion
of Gishboliner, Glock, Michaeli and Sgueglia~\cite[Section~4]{GGMS}.
Our proofs build on the convex compression, integer-threshold
rigidity, profile switching and exceptional-vertex completion
mechanisms of Liu, Turevskii, Wang and Yan~\cite{LTWY}; these
mechanisms are not claimed as new here. Weak hypergraph regularity
has its origins in the work of Chung~\cite{Chung}; its use for
embedding linear hypergraphs was developed by Kohayakawa,
Nagle, R\"odl and Schacht~\cite{KNRS}. Fractional matching and
covering are standard linear programming tools; we refer to
Schrijver~\cite{Schrijver} for a systematic treatment. We give proofs of the
rigidity and completion statements in the forms needed for
our results. Section~\ref{sec:clustering} proves the covering-weight
lemma, Section~\ref{sec:stability-proof} proves stability, and
Section~\ref{sec:exact-proof} derives the exact normal form.

The defect and random-host extensions form another direction.
Gishboliner, Glock, Michaeli and Sgueglia~\cite{GGMS} establish
an asymptotic AFL bound for almost complete hypergraphs and
transfer it to sparse random hypergraphs. In the graph setting,
Gishboliner, Krivelevich and Michaeli~\cite{GKM} obtained a
random-host matching result in their study of colour-biased
Hamilton cycles. These results concern the persistence of large
monochromatic matchings when the host is incomplete; our
theorems describe near-extremal and exact colourings of the
complete host.

For graphs, Xu, Yang and Zhang~\cite[Theorem~2]{XYZ} characterized
the critical complete-graph colourings for matching Ramsey
numbers, allowing unequal prescribed matching sizes. Their proof
uses the Gallai--Edmonds decomposition; see Lov\'asz and
Plummer~\cite{LP} for the underlying graph matching theory.
Our exact statement concerns
all sufficiently large orders for fixed $r,q$. No explicit value
of $n_{\mathrm{ex}}(r,q)$ is supplied.

We write $\one=(1,\ldots,1)$ and $\ee_i$ for the $i$th coordinate
vector. The support of a probability measure is its closed
topological support. All asymptotic statements keep $r,q$ fixed. A partition may have
empty parts unless nonemptiness is explicitly required. For a
partition $A_1,\ldots,A_q$ and an edge $e$, the vector
$(|e\cap A_1|,\ldots,|e\cap A_q|)$ is the \emph{profile} of $e$.

\section{Rigidity of near-critical covering weights}\label{sec:clustering}

The following finite statement is the only rigidity input used
in the stability proof. It says that near-optimal covering weights concentrate near
a finite collection of axis points with integer reciprocals.
The compactness argument below proves the required uniform
finite statement without providing an explicit bound on its parameters.

\begin{lemma}[Robust cover clustering]\label{lem:clustering}
For every $\alpha>0$ there is $\delta>0$ such that the following
holds. Let $t$ be a positive integer and let
$z_1,\ldots,z_t\in[0,1]^q$ satisfy
\begin{equation}\label{eq:array-means}
 \frac1t\sum_{v=1}^t z_v^{(i)}\le\frac1L+\delta
 \qquad(i\in[q]).
\end{equation}
For independent uniform indices $v_1,\ldots,v_r\in[t]$, sampled
with replacement, suppose that
\begin{equation}\label{eq:array-failure}
 \Prob\left(\max_i\sum_{a=1}^r z_{v_a}^{(i)}<1\right)\le\delta.
\end{equation}
Put $\kappa=1/(4r^2)$. There is $c\in\cC_{r,q}$ such that the
sets
$ I_i=\{v:\norminf{z_v-\ee_i/c_i}<\kappa\}\ (i=1,\dots,q)$ and
 $I_0=[t]\setminus\bigcup_{i=1}^q I_i$
form a partition of $[t]$ and satisfy
\begin{equation}\label{eq:array-clusters}
 \left|\frac{|I_i|}{t}-\frac{c_i}{L}\right|\le\alpha
 \quad(i\in[q]),\qquad |I_0|\le\alpha t.
\end{equation}
\end{lemma}

\begin{proof}
Suppose the lemma fails for some $\alpha>0$. For each positive
integer $h$, take a counterexample with $\delta=1/h$ and form
its empirical measure
$\sigma_h=t_h^{-1}\sum_{v=1}^{t_h}\delta_{z_{h,v}}$, where
$\delta_z$ denotes the unit point mass at $z$.
Probability measures on the compact cube have a weakly
convergent subsequence; this is the compact-space case of
Prokhorov's theorem, as presented in Billingsley~\cite{Billingsley}. Write $\sigma_h\Rightarrow\sigma$ along
such a subsequence and $S=\supp\sigma$.

\smallskip\noindent
\emph{Step 1: an exact covering problem on the support.}
Coordinate functions are continuous, so~\eqref{eq:array-means}
implies
\begin{equation}\label{eq:limit-means}
 \int z_i\,d\sigma(z)\le1/L\qquad(i\in[q]).
\end{equation}
We claim that, for all $z^{(1)},\ldots,z^{(r)}\in S$,
including repetitions,
\begin{equation}\label{eq:support-cover}
 \max_i\sum_{a=1}^r z_i^{(a)}\ge1.
\end{equation}
If this failed, continuity would give relative open
neighbourhoods $U_a$ of $z^{(a)}$ whose entire Cartesian product
violates the inequality. Each $\sigma(U_a)>0$ because
$z^{(a)}\in S$. The Portmanteau theorem~\cite{Billingsley} gives
$\liminf_h\sigma_h(U_a)\ge\sigma(U_a)>0$. Independent sampling
would therefore violate the inequality with probability at
least $\prod_a\sigma_h(U_a)$, bounded away from zero.
This contradicts~\eqref{eq:array-failure}. The same argument
works when some of the support points coincide.

\smallskip\noindent
\emph{Step 2: the covering condition forces integer thresholds.}
Let $C=\conv S$. Since $S$ is compact and the ambient space
is finite-dimensional, $C$ is compact. Choose $u\in C$ minimizing
$\max_i u_i$. The barycentre $\int z\,d\sigma(z)$ belongs to $C$,
since $C$ is closed and convex and contains the support of $\sigma$. Thus
for each $i\in[q]$, we have
 $u_i\le1/L$.
 
If $C$ has affine dimension at most $q-1$, Carath\'eodory's
theorem~\cite{Rockafellar} represents $u$ using at most $q$
points of $S$.
If $C$ has dimension $q$, then $u$ is on its boundary:
otherwise $u-\varepsilon\one\in C$ for small $\varepsilon>0$,
contrary to minimality. Take a supporting hyperplane through
$u$, using the supporting-hyperplane theorem~\cite{Rockafellar}. In any finite convex representation of $u$ by points of $S$,
every point with a positive coefficient lies in that hyperplane,
because all points of $C$ lie on one side.
Carath\'eodory's theorem in this face again uses at most $q$
points. Thus, omitting zero coefficients, we have
$ u=\sum_{j=1}^p\lambda_j w_j,\ 
p\le q,\ w_j\in S,\ 
 \lambda_j>0$, and $\sum_j\lambda_j=1$.

Set $b_j=\ceil{L\lambda_j}-1$, a nonnegative integer strictly
smaller than $L\lambda_j$. If $\sum_j b_j\ge r$, choose
integers $0\le a_j\le b_j$ with $\sum_j a_j=r$.
Apply~\eqref{eq:support-cover} to $a_j$ copies of each $w_j$.
For some coordinate $\ell$,
$ 1\le\sum_j a_jw_{j,\ell}
   <\sum_j L\lambda_jw_{j,\ell}
   =Lu_\ell\le1$,
a contradiction. The middle inequality is strict because every
$a_j<L\lambda_j$, all coordinates are nonnegative, and the
first sum has a positive term. Therefore
$ r-1\ge\sum_j b_j
 =\sum_j\ceil{L\lambda_j}-p
 \ge L-p\ge L-q=r-1$.
Equality holds throughout. Hence $p=q$ and
$\sum_j\ceil{L\lambda_j}=\sum_j L\lambda_j$, so every
$c_j:=L\lambda_j$ is a positive integer. Their sum is $L$.

\smallskip\noindent
\emph{Step 3: identify the support and its masses.}
For each $j$, the profile $c-\one+\ee_j$ is nonnegative and
has sum $r$. By~\eqref{eq:support-cover}, some coordinate
$\ell_j$ satisfies
$ 1\le Lu_{\ell_j}-\sum_{a\ne j}w_{a,\ell_j}\le1$.
Equality throughout gives $Lu_{\ell_j}=1$ and
$w_{a,\ell_j}=0$ for $a\ne j$. Since
$Lu_{\ell_j}=\sum_a c_aw_{a,\ell_j}$, it also gives
$c_jw_{j,\ell_j}=1$. The coordinates $\ell_j$ are all distinct:
if $\ell_j=\ell_a$ with $j\ne a$, the latter equality for $j$
contradicts the zero forced by the equality for $a$.
Since $j\mapsto\ell_j$ is a permutation of $[q]$, relabel
the selected support points and their associated integers so
that $\ell_j=j$. All $q$ coordinates have been used, and therefore for every $j\in[q]$, we have
$ w_j=\ee_j/c_j$.

For any $z\in S$, apply~\eqref{eq:support-cover} to the following
$r$ support points: $z$ and
$c_i-1$ copies of $\ee_i/c_i$ for each $i$. In coordinate $i$, their sum is $z_i+(c_i-1)/c_i$.
Thus $c_i z_i\ge1$ for some $i$. Consequently
\[
 g(z):=\sum_i c_i z_i\ge1\quad(z\in S),\qquad
 \int g\,d\sigma\le\sum_i c_i/L=1
\]
by~\eqref{eq:limit-means}. Thus $g=1$ on $S$: a support point
with $g>1$ would, by continuity, have a positive-measure
neighbourhood on which $g$ exceeds one by a fixed amount.
Combining $g(z)=1$ with $c_i z_i\ge1$ and nonnegativity gives
$z=\ee_i/c_i$. If $p_i$ is the mass at $\ee_i/c_i$, then
$p_i/c_i\le1/L$ by~\eqref{eq:limit-means}. Since
$\sum_i p_i=1=\sum_i c_i/L$, all these bounds are equalities:
\begin{equation}\label{eq:axis-law}
 \sigma=\sum_i\frac{c_i}{L}\delta_{\ee_i/c_i}.
\end{equation}

Finally, the open balls of radius $\kappa$ around the axis points
are disjoint: their centres have $\ell_\infty$-distance at least
$1/r>2\kappa$. Their boundaries contain no atom
of~\eqref{eq:axis-law}. Weak convergence gives
$|I_i|/t_h\to c_i/L$ and $|I_0|/t_h\to0$.
For large $h$ this is~\eqref{eq:array-clusters}, contradicting
the choice of the counterexamples.
\end{proof}

\section{Proof of stability}\label{sec:stability-proof}

We first obtain covering weights from a regular partition.
Lemma~\ref{lem:clustering} then identifies the parts, and an
elementary repair puts their sizes within the exact capacities.

\subsection{Regularity and fractional matching lift}
For disjoint nonempty sets $X_1,\ldots,X_r$, let
$d_H(X_1,\ldots,X_r)$ be the number of edges of $H$ meeting
each set once, divided by $|X_1|\cdots|X_r|$.
The tuple is \emph{$\eta$-regular in $H$} if
\[
 |d_H(Y_1,\ldots,Y_r)-d_H(X_1,\ldots,X_r)|\le\eta
\]
whenever $Y_a\subseteq X_a$ and $|Y_a|\ge\eta|X_a|$ for all $a$.
A partition is equitable if its part sizes differ by at most one.

This is the vertex-partition notion of weak hypergraph regularity
studied in~\cite{Chung,KNRS}. We use the following dense
multicolour form. It is~\cite[Lemma~3.1]{GGMS} with $p=1$, where the upper
uniformity assumptions hold automatically with, for example,
upper-density bound $2$. If needed, increase the prescribed
lower bound on the number of parts before applying that lemma.

\begin{lemma}[Weak regularity]\label{lem:regularity}
For every $\eta>0$ and integer $t_0\ge r$, there are integers
$T,N$ such that any $r$-uniform hypergraphs $H_1,\ldots,H_q$
on a common vertex set of size $n\ge N$ admit an equitable
partition $X_1,\ldots,X_t$, with $t_0\le t\le T$, such that,
for every $i\in[q]$, at most $\eta\binom tr$ of the unordered
$r$-tuples of distinct parts fail to be $\eta$-regular in $H_i$.
\end{lemma}

Call a tuple \emph{jointly regular} if it is regular in every
colour. At most $q\eta\binom tr$ tuples are not jointly regular.
Fix $0<\eta<d<1/q$. Given an equitable partition of a colouring,
define an $r$-graph $J_i$ on $[t]$ for each $i\in[q]$ by
\begin{equation}\label{eq:reduced}
 Q\in E(J_i)
 \ \Longleftrightarrow\
 \begin{cases}
  Q\in\binom{[t]}r\text{ is jointly regular},\\
  d_{H_i}((X_v)_{v\in Q})\ge d.
 \end{cases}
\end{equation}
All colours of density at least $d$ are retained; the $J_i$ may
overlap. This is important: an edge of density at least $d$
in colour $i$ must be retained in $J_i$, even if it also belongs
to another $J_j$. Since the colour densities sum to one, every
jointly regular tuple belongs to at least one $J_i$.

A \emph{fractional matching} in an $r$-graph $J$ assigns
nonnegative weights to its edges so that the sum at each vertex
is at most one. A \emph{fractional vertex cover} assigns
nonnegative weights to its vertices so that the sum on every
edge is at least one. Linear programming duality, in its standard
fractional packing--covering form~\cite{Schrijver}, gives
\begin{align}
 \nu^*(J)
 &=\max\left\{\sum_{Q\in E(J)}x_Q:
       x_Q\ge0,\ \sum_{Q\ni v}x_Q\le1\text{ for all }v\right\}
       \notag\\
 &=\min\left\{\sum_{v\in V(J)}z_v:
       0\le z_v\le1,\ \sum_{v\in Q}z_v\ge1
                           \text{ for all }Q\in E(J)\right\}.
 \label{eq:duality}
\end{align}
The upper bounds $z_v\le1$ cause no loss: truncating a larger
coordinate at one preserves every cover constraint.

\begin{lemma}[Fractional matching lift]\label{lem:lift}
Let $J_i$ be defined by~\eqref{eq:reduced}, with
$0<\eta<d<1/q$. Suppose the clusters have size $b$ or $b+1$,
where $b$ is a positive integer. Then
\begin{equation*}\label{eq:lift}
 \nu(H_i)\ge(1-2\eta)b\,\nu^*(J_i)-\binom tr
 \qquad(i\in[q]).
\end{equation*}
\end{lemma}

\begin{proof}
Take an optimal fractional matching $(x_Q)$ in $J_i$. For each
$Q$, set a demand of $N_Q=\floor{(1-2\eta)b x_Q}$ colour-$i$
edges across its clusters. Process these demands in any order,
choosing every edge disjoint from all edges already chosen.
The total demand on a cluster $X_v$ is at most
\[
 \sum_{Q\ni v}N_Q\le(1-2\eta)b\sum_{Q\ni v}x_Q
                         \le(1-2\eta)b.
\]
At every stage each cluster has at least $2\eta b$ unused
vertices. Since $b\ge1$ gives $2\eta b\ge\eta(b+1)$, each unused set is
large enough for regularity. For the current tuple the
colour-$i$ density on its unused sets is at least $d-\eta>0$,
so the next edge exists. The resulting matching has size
\[
 \sum_Q N_Q\ge(1-2\eta)b\sum_Qx_Q-|E(J_i)|
             \ge(1-2\eta)b\,\nu^*(J_i)-\binom tr.
\]
The total rounding loss is therefore at most the number of
edges of the reduced hypergraph.
\end{proof}

\subsection{Repairing the part sizes}
\begin{lemma}[Partition repair]\label{lem:repair}
Let $c\in\cC_{r,q}$ and $V=A_1\sqcup\cdots\sqcup A_q$.
Suppose at most $B$ edges of $\chi$ violate
$\chi(e)=i\Longrightarrow |e\cap A_i|\ge c_i$.
Put
$ E=\sum_i\max\{0,|A_i|-(c_i m-1)\}$.
Then
\begin{equation}\label{eq:repair}
 \dist(\chi,\cE_{r,q}(V))\le B+E\binom{n-1}{r-1}.
\end{equation}
\end{lemma}

\begin{proof}
The capacities sum to $n+D$, so the total unused capacity in
underfull parts is $E+D\ge E$. Move exactly $E$ vertices from
overfull parts to underfull parts without exceeding a receiving
capacity. Edges avoiding the moved vertices retain their
intersection profiles. Thus at most the original $B$ violations
and the edges meeting a moved vertex need recolouring.
Recolour every such edge whose original colour is ineligible
with any eligible colour. Such a colour exists
by~\eqref{eq:threshold-range}, and the new partition respects
the capacities. At most $E\binom{n-1}{r-1}$ edges meet a
moved vertex, proving~\eqref{eq:repair}.
\end{proof}

\subsection{Proof of the stability theorem}
\begin{proof}[Proof of Theorem~\ref{thm:stability}]
We may assume $0<\eps\le1$. Throughout this proof, the
parameters in~\eqref{eq:parameters} are computed from the current order $n$.

\smallskip\noindent
\emph{Step 1: parameter selection.}
Set $\alpha=\eps/[4(4q+3)]$ and obtain $\delta>0$ from
Lemma~\ref{lem:clustering}; decrease it so that $\delta\le1/L$.
Choose
$ 0<d\le\min\left\{\frac1{2q},\frac{\eps}{2^{r+3}q}\right\}$.
Next choose $\eta,\mu>0$ so that
$\eta<d,\ q\eta\le\delta/2$, and
 $\frac{1+\mu}{L(1-2\eta)}<\frac1L+\frac{\delta}{2}$.
This is possible because the last fraction tends to $1/L$
as $\eta,\mu\to0$. Take an integer
$ t_0\ge\max\left\{r,\frac8\eps,\frac{2\binom r2}{\delta}\right\}$,
and apply Lemma~\ref{lem:regularity} to obtain $T,N$.
Finally choose $n_0\ge N$ so large that, for every $n\ge n_0$,
$ T\le\alpha n,\ q\le\alpha n$, and
$ \frac{(1+\mu)/L+T^r/n}{(1-2\eta)(1-T/n)}
                 \le\frac1L+\delta$.
The last requirement is possible because its left side tends
to $(1+\mu)/(L(1-2\eta))$, which is less than $1/L+\delta/2$.
Thus all constants, including $T$, have been fixed before $n_0$.

\smallskip\noindent
\emph{Step 2: obtain the threshold vector from cover weights.}
Let $\chi$ satisfy the theorem's assumptions, and take a regular
partition with $t_0\le t\le T$ and cluster sizes $b$ or $b+1$,
where $b=\floor{n/t}$. Define all the $J_i$ by~\eqref{eq:reduced}
and choose optimal covers $z_v^{(i)}\in[0,1]$.
By~\eqref{eq:duality} and Lemma~\ref{lem:lift},
\begin{align*}
 \frac1t\sum_v z_v^{(i)}
 =\frac{\nu^*(J_i)}t
 \le\frac{(1+\mu)n/L+\binom tr}{(1-2\eta)bt}
 \le\frac{(1+\mu)/L+T^r/n}{(1-2\eta)(1-T/n)}
 \le\frac1L+\delta.
 \label{eq:cover-means}
\end{align*}
Here $bt\ge n-t\ge n-T$.
For independent uniform indices $v_1,\ldots,v_r$, the probability
of a repetition is at most $\binom r2/t$. Conditional on being
distinct, their underlying $r$-set is uniform. If it is jointly
regular, it belongs to some $J_i$ and hence has cover sum at
least one in that coordinate. Thus, for
$z_v=(z_v^{(1)},\ldots,z_v^{(q)})$,
$ \Prob\left(\max_i\sum_a z_{v_a}^{(i)}<1\right)
 \le q\eta+\frac{\binom r2}{t_0}\le\delta$.
Lemma~\ref{lem:clustering} now gives $c\in\cC_{r,q}$ and
$I_0,I_1,\ldots,I_q$ satisfying~\eqref{eq:array-clusters}.

Let $A_i^0=\bigcup_{v\in I_i}X_v$ and
$Z=\bigcup_{v\in I_0}X_v$. Equitability gives
$ |Z|\le\alpha t(b+1)\le2\alpha n$, and
$ \left||A_i^0|-\frac{c_i n}{L}\right|
 \le\alpha n+t\le2\alpha n$.
For the second bound, replacing each cluster size by $n/t$
changes the sum over $I_i$ by at most $t$.
Assign $Z$ arbitrarily to the parts, obtaining a partition
$A_1,\ldots,A_q$ of $V$ with
\begin{equation}\label{eq:approximate-parts}
 \left||A_i|-\frac{c_i n}{L}\right|\le4\alpha n.
\end{equation}

\smallskip\noindent
\emph{Step 3: bound every violation of the threshold rule.}
Consider a jointly regular tuple $Q$ of clusters in
$I_1\cup\cdots\cup I_q$, and put $a_i=|Q\cap I_i|$.
By the definition of the clusters $I_j$, each $z_v^{(i)}$
differs from $1/c_i$ by less than $\kappa$ when $v\in I_i$,
and is less than $\kappa$ otherwise. Hence
$ \sum_{v\in Q}z_v^{(i)}\le \frac{a_i}{c_i}+r\kappa$.
If $a_i<c_i$, then, since $c_i\le r$,
$ \sum_{v\in Q}z_v^{(i)}
 \le1-\frac1{c_i}+r\kappa
 \le1-\frac1r+\frac1{4r}<1$.
The cover constraint rules out $Q\in J_i$. Joint regularity
and~\eqref{eq:reduced} therefore imply that its colour-$i$
density is less than $d$. Summing over ineligible colours,
the wrong-colour density on this tuple is at most $qd$.

Every threshold violation either meets $Z$, uses a cluster twice,
lies in a tuple that is not jointly regular, or is one of the
low-density wrong-colour edges just counted. Their total number
$B$ consequently satisfies
\begin{align*}
 B\le{}
 |Z|\binom{n-1}{r-1}
 +\sum_{v=1}^t\binom{|X_v|}{2}\binom{n-2}{r-2}
 +q(\eta+d)\binom tr(b+1)^r.
 \label{eq:all-errors}
\end{align*}
Here $b+1\le2n/t$ and
$\sum_v\binom{|X_v|}{2}\le t(b+1)^2/2\le2n^2/t$.
Thus the middle term is at most $2n^r/t$,
and the last term is at most $2^r q(\eta+d)n^r$. Therefore
$ \frac{B}{n^r}\le2\alpha+\frac2{t_0}+2^r q(\eta+d)$.

Finally, $m=(n+q+D)/L\ge n/L$, so $u_i=c_i m-1\ge c_i n/L-1$.
Equation~\eqref{eq:approximate-parts} gives
\[
 E=\sum_i\max\{0,|A_i|-u_i\}
 \le q(4\alpha n+1)\le(4q+1)\alpha n.
\]
Applying Lemma~\ref{lem:repair} and
$\binom{n-1}{r-1}\le n^{r-1}$, we obtain
\[
 \frac{\dist(\chi,\cE_{r,q}(V))}{n^r}
 \le(4q+3)\alpha+\frac2{t_0}+2^r q(\eta+d)
 \le\frac{3\eps}{4}<\eps.
\]
The three terms are each at most $\eps/4$ by the choices of
$\alpha,t_0,d$ and $\eta<d$. This proves the theorem.
\end{proof}

\section{From stability to the exact normal form}\label{sec:exact-proof}

We prove sufficiency first. For necessity, we isolate one
completion inequality: its empty-matching case gives the part
capacities, and its nonempty cases give the deficit budgets.

\subsection{Sufficiency and bounded witnesses}
\begin{proof}[Proof of the sufficiency in Theorem~\ref{thm:exact}]
Suppose~\eqref{eq:normal-sizes}--\eqref{eq:budget} hold.
If $N$ were a colour-$i$ matching with $m$ edges, then
\[
 |W_i|
 \ge\sum_{e\in N}|e\cap W_i|
 \ge c_i m-\sum_{e\in N\cap\mathcal F_i}\delta_i(e)
 \ge c_i m-s_i=|W_i|+1,
\]
a contradiction. Thus every colour has matching number at most
$m-1=k$, and~\eqref{eq:AFL} makes the colouring extremal.

To justify the bounded tests, take an inclusion-minimal matching
violating~\eqref{eq:budget}. If it had more than $s_i+1$ edges,
any $s_i+1$ of them would already have total deficit at least
$s_i+1$, since the deficits are positive integers. This would
contradict minimality. Thus a violation has a witness of at most
$s_i+1\le D+1\le L$ edges.
\end{proof}

\subsection{A profile switching lemma}
Let $a=(a_1,\ldots,a_q)$ be a nonnegative integer vector of sum
$r$, and let $h$ be a positive integer. For disjoint sets $U_j$
of sizes $a_jh$, write
\[
 \mathcal K(a;U_1,\ldots,U_q)
 =\{e\subseteq\textstyle\bigcup_j U_j:|e\cap U_j|=a_j
                                               \text{ for all }j\}.
\]
For an edge family $\mathcal B$, let
$d_{\mathcal B}(v)=|\{e\in\mathcal B:v\in e\}|$ and
$\Delta(\mathcal B)=\max_vd_{\mathcal B}(v)$.
The following form of the switching argument
of~\cite{LTWY} gives an explicit sufficient bound on the
maximum degree of a forbidden edge family. In the proof,
we call the edges of $\mathcal B$ \emph{bad} and the other
edges \emph{good}.

\begin{lemma}[Profile switching]\label{lem:switching}
Allow $U_j=\varnothing$ when $a_j=0$. If $h\ge r$ and $\mathcal B\subseteq\mathcal K(a;U_1,\ldots,U_q)$
satisfies
\begin{equation}\label{eq:switching-bound}
 r\Delta(\mathcal B)<(h-1)_{r-1}
 :=(h-1)(h-2)\cdots(h-r+1),
\end{equation}
then $\mathcal K(a;U_1,\ldots,U_q)\setminus\mathcal B$ contains
a perfect matching of size $h$. In particular,
$\Delta(\mathcal B)=o(h^{r-1})$ suffices as $h\to\infty$.
\end{lemma}

\begin{proof}
For each $j$, write $U_j$ as the disjoint union of $h$
labelled blocks of size $a_j$ (all empty if $a_j=0$).
Unite blocks with the same label. This gives a perfect
matching in the complete profile hypergraph. Choose one,
$\mathcal M$, with as few bad edges as possible, and suppose
$e_0\in\mathcal M\cap\mathcal B$.

Divide the positions $0,\ldots,r-1$ into groups of sizes
$a_1,\ldots,a_q$. Label the vertices of every edge of
$\mathcal M$ by these positions, so that the same positions
always belong to the same part $U_j$. Fix all labels before
making the following choices.

Choose ordered distinct edges
$e_1,\ldots,e_{r-1}\in\mathcal M\setminus\{e_0\}$.
For $0\le a,b<r$, let $x_{a,b}$ be the vertex in position $b$
of $e_a$. Put
\[
 f_s=\{x_{a,b}:a+b\equiv s\pmod r\},\qquad 0\le s<r.
\]
Each new edge has one vertex from each position, hence the
correct profile. The new edges are disjoint, cover the same
vertices as the selected old edges, and $f_s$ contains the
fixed vertex $x_{0,s}$.

For a fixed $s$, any prescribed edge arising as $f_s$ determines
at most one ordered tuple $(e_1,\ldots,e_{r-1})$ of old edges. Each of its
vertices identifies its old matching edge and its fixed
position $b$, and hence that edge's row $a\equiv s-b\pmod r$.
At most $d_{\mathcal B}(x_{0,s})$ choices therefore make $f_s$
bad. At most $r\Delta(\mathcal B)$ choices make any new edge
bad, whereas there are $(h-1)_{r-1}$ choices in total.
By~\eqref{eq:switching-bound}, one choice has all new edges
good. Switching to these edges reduces the number of bad
edges, a contradiction.
\end{proof}

\subsection{One inequality for capacities and deficits}
For a matching $M$, write $V(M)=\bigcup_{e\in M}e$, with
$V(\varnothing)=\varnothing$. The next lemma refines the
exceptional-vertex completion in~\cite{LTWY} by allowing an arbitrary bounded matching to be prescribed first.
Its contrapositive is the main idea: if a prescribed matching
uses too few vertices of $W$, then the remaining vertices of $W$
allow it to be extended to a matching of size $m$.

\begin{lemma}[Bounded matching completion]\label{lem:completion}
Fix $c\in\cC_{r,q}$, $i\in[q]$, an integer $B\ge0$, and
$\gamma>0$. Consider a sequence of $q$-colourings of complete $r$-graphs,
with an integer parameter
$m\to\infty$ and vertex partitions
\[
 V=A_1\sqcup\cdots\sqcup A_q\sqcup X,\qquad
 |A_j|=c_jm+o(m),\qquad |X|=o(m).
\]
Let $p^{(i)}=c-\one+\ee_i$. Suppose the family of edges in
the cores $A_1,\ldots,A_q$ with profile $p^{(i)}$ that are
not coloured $i$ has maximum degree $o(m^{r-1})$.
Let $Y\subseteq X$ and put $W=A_i\sqcup Y$.
Assume each $x\in Y$ belongs to at least $\gamma m^{r-1}$
colour-$i$ edges $\{x\}\cup S$ with
\begin{equation}\label{eq:baseline}
 S\subseteq A_1\cup\cdots\cup A_q,\qquad
 |S\cap A_j|=c_j-1\quad(j\in[q]).
\end{equation}
Suppose $\nu(H_i)<m$ for every member of the sequence. Then,
for all sufficiently large members
of the sequence, every colour-$i$ matching $M$ with $|M|\le B$,
including the empty matching, satisfies
\begin{equation}\label{eq:completion}
 c_i|M|-|V(M)\cap W|\le c_i m-|W|-1.
\end{equation}
\end{lemma}

The left side of~\eqref{eq:completion} may be negative for a
general matching. It equals the sum of the positive deficits
when every edge of $M$ is violating relative to $W$.

\begin{proof}
Fix a member of the sequence and a colour-$i$ matching $M$
with $|M|\le B$. Set
\[
 \ell=|M|,\quad a_j=|V(M)\cap A_j|,\quad
 b=|V(M)\cap Y|,\quad t=|Y|-b.
\]
Then $\ell\le B$, $a_j\le rB$ and $t=o(m)$. We work with
a sufficiently large member of the sequence; the estimates
below will be uniform over all such $M$.
If~\eqref{eq:completion} fails, integrality gives
\begin{equation}\label{eq:completion-failure}
 c_i\ell-a_i-b\ge c_i m-|W|.
\end{equation}

For each $x\in Y\setminus V(M)$, choose a colour-$i$ edge
$\{x\}\cup S_x$ satisfying~\eqref{eq:baseline}, disjoint from
$M$ and the previously chosen edges. This is possible greedily.
At most $rB+(r-1)|Y|=o(m)$ core vertices are forbidden.
Each forbidden core vertex excludes at most
$\binom{|V|-2}{r-2}=O_{r,q}(m^{r-2})$ candidates through $x$,
since $|V|=Lm+o(m)$. Thus at most $o(m^{r-1})$ of the at least
$\gamma m^{r-1}$ colour-$i$ candidates are excluded.
The noncore vertex of each edge is its own distinct $x$,
outside $V(M)$.

After these $t$ choices, we seek
$h=m-\ell-t=(1-o(1))m$ further colour-$i$ edges of profile
$p^{(i)}$. In the $i$th core, the number of remaining vertices
minus the number required is
\begin{align*}
 |A_i|-a_i-t(c_i-1)-c_i h
 &=c_i\ell-a_i+t+|A_i|-c_i m\notag\\
 &=c_i\ell-a_i-b+|W|-c_i m\ge0,
 \label{eq:own-surplus}
\end{align*}
by~\eqref{eq:completion-failure}. For $j\ne i$, the difference is
\begin{align*}
 |A_j|-a_j-t(c_j-1)-(c_j-1)h
 &=|A_j|-(c_j-1)m-a_j+(c_j-1)\ell\notag\\
 &=m+o(m)>0.
\end{align*}
For sufficiently large $m$, we have $h\ge r$. Choose subsets
of the remaining cores having exactly
$p_j^{(i)}h$ vertices in each $A_j$. The wrong-colour family
on them still has maximum degree
$o(m^{r-1})=o(h^{r-1})$. Lemma~\ref{lem:switching} supplies
$h$ disjoint colour-$i$ edges. Together with $M$ and the
$t$ selected links, they form a matching of size
$\ell+t+h=m$, a contradiction.
The number of forbidden core vertices is bounded by
$rB+(r-1)|Y|$, independently of the choice of $M$. The switching
bound also depends only on the given wrong-colour maximum degree
and on $h=m-O(B)-o(m)$. Consequently, one sufficiently large
index works simultaneously for all matchings with $|M|\le B$.
\end{proof}

\subsection{Necessity of the normal form}
\begin{proof}[Proof of the necessity in Theorem~\ref{thm:exact}]
Suppose that there are counterexamples for arbitrarily large
$n$. Choose such a sequence with $n\to\infty$, and use
$k,m,D$ from~\eqref{eq:parameters}. All asymptotic estimates
below refer to this sequence.

\smallskip\noindent
\emph{Step 1: remove vertices with large error degree.}
Extremality gives $\max_i\nu(H_i)\le k=n/L+O_{r,q}(1)$.
The sequential form of Theorem~\ref{thm:stability} yields
$\psi\in\cE_{r,q}(V)$ with
\[
 \mathcal B=\{e:\chi(e)\ne\psi(e)\},\qquad
 |\mathcal B|=o(n^r)=o(m^r).
\]
Let $A_1',\ldots,A_q'$ be a witnessing partition for $\psi$.
After passing to a subsequence, its vector $c\in\cC_{r,q}$ is
fixed, because there are finitely many such vectors.
By~\eqref{eq:part-sizes}, $|A_j'|=c_jm+O_{r,q}(1)$.

Write $\beta=|\mathcal B|/m^r=o(1)$, put
$\theta=\sqrt\beta+m^{-1/2}$, and let
$ X=\{v:d_{\mathcal B}(v)>\theta m^{r-1}\}$ and
 $A_j=A_j'\setminus X$.
Since $\sum_v d_{\mathcal B}(v)=r|\mathcal B|$, we have
$ |X|\le\frac{r\beta}{\theta}m=o(m)$ and
 $|A_j|=c_jm+o(m)$.
An edge in the cores with profile $p^{(i)}=c-\one+\ee_i$
has only colour $i$ eligible in $\psi$. Every such edge not
coloured $i$ in $\chi$ therefore belongs to $\mathcal B$.
Its wrong-colour family has maximum degree at most
$\theta m^{r-1}=o(m^{r-1})$.

\smallskip\noindent
\emph{Step 2: assign the exceptional vertices to colours.}
For $x\in X$, the number of edges $\{x\}\cup S$ with
$|S\cap A_j|=c_j-1$ for every $j$ is
\[
 \prod_{j=1}^q\binom{|A_j|}{c_j-1}
       =\Theta_{r,q}(m^{r-1}).
\]
Here $|A_j|=c_jm+o(m)$ and $\sum_j(c_j-1)=r-1$.
By the pigeonhole principle, one colour $\tau(x)$ occurs on at least
$\gamma m^{r-1}$ of these edges, for a fixed $\gamma>0$
depending only on $r,q$. Indeed, the leading coefficient of
the product above is $\prod_j c_j^{c_j-1}/(c_j-1)!>0$, with
$0!=1$; there are only finitely many possible vectors $c$.
Fix such a colour for each $x$ and put
\[
 X_i=\{x\in X:\tau(x)=i\},\qquad W_i=A_i\sqcup X_i.
\]
These sets form a partition of $V$.

\smallskip\noindent
\emph{Step 3: recover the capacities and deficit budgets.}
Apply Lemma~\ref{lem:completion} with $B=L$, $Y=X_i$ and
$W=W_i$, for each $i$. Its hypotheses follow from Steps~1--2
and $\nu(H_i)\le m-1$. Since there are only $q$ colours, the
conclusions hold simultaneously for all $i$ from some index
onwards. Hence every colour-$i$ matching with
at most $L$ edges obeys
\begin{equation}\label{eq:small-budget}
 c_i|M|-|V(M)\cap W_i|\le c_i m-|W_i|-1.
\end{equation}
Taking $M=\varnothing$ defines nonnegative integers
$ s_i=c_i m-|W_i|-1\ge0$.
Because the $W_i$ partition $V$,
\begin{equation*}\label{eq:exact-size-recovery}
 |W_i|=c_i m-1-s_i,\qquad
 \sum_i s_i=Lm-n-q=D.
\end{equation*}
For a matching $M$ of violating colour-$i$ edges, disjointness gives
$ \sum_{e\in M}\delta_i(e)=c_i|M|-|V(M)\cap W_i|$.
If it violated the budget $s_i$, the bounded-witness argument
at the start of this section would give a violating submatching
with at most $s_i+1\le D+1\le L$ edges. This
contradicts~\eqref{eq:small-budget}.

Thus~\eqref{eq:normal-sizes}--\eqref{eq:budget} hold for every
sufficiently large member of the counterexample sequence.
This contradiction proves the existence of
$n_{\mathrm{ex}}(r,q)$ and the necessity of the normal form.
\end{proof}

\section*{Usage of AI tools}
The proofs presented in this paper were developed by ChatGPT.
The authors take full responsibility for the correctness of
the mathematical arguments and for the entire content of
the manuscript.


\begin{thebibliography}{99}
	
	\bibitem{AK}
	R. Ahlswede and L. H. Khachatrian,
	The complete intersection theorem for systems of finite sets,
	\emph{European J. Combin.} \textbf{18} (1997), 125--136.
	
	\bibitem{AFL}
	N. Alon, P. Frankl and L. Lov\'asz,
	The chromatic number of Kneser hypergraphs,
	\emph{Trans. Amer. Math. Soc.} \textbf{298} (1986), 359--370.
	
	\bibitem{Barany78}
	I. B\'ar\'any,
	A short proof of Kneser's conjecture,
	\emph{J. Combin. Theory Ser. A} \textbf{25} (1978), 325--326.
	
	\bibitem{Billingsley}
	P. Billingsley,
	\emph{Convergence of Probability Measures}, 2nd ed.,
	Wiley Series in Probability and Statistics,
	John Wiley \& Sons, New York, 1999.
	
	\bibitem{Chung}
	F. R. K. Chung,
	Regularity lemmas for hypergraphs and quasi-randomness,
	\emph{Random Structures Algorithms} \textbf{2} (1991), 241--252.
	
	\bibitem{CL}
	E. J. Cockayne and P. J. Lorimer,
	The Ramsey number for stripes,
	\emph{J. Austral. Math. Soc.} \textbf{19} (1975), 252--256.
	
	\bibitem{Erdos65}
	P. Erd\H{o}s,
	A problem on independent $r$-tuples,
	\emph{Ann. Univ. Sci. Budapest. E\"otv\"os Sect. Math.}
	\textbf{8} (1965), 93--95.
	
	\bibitem{Frankl13}
	P. Frankl,
	Improved bounds for Erd\H{o}s' matching conjecture,
	\emph{J. Combin. Theory Ser. A} \textbf{120} (2013), 1068--1072.
	
	\bibitem{FK21}
	P. Frankl and A. Kupavskii,
	Beyond the Erd\H{o}s matching conjecture,
	\emph{European J. Combin.} \textbf{95} (2021), Paper No. 103338.
	
	\bibitem{FK22}
	P. Frankl and A. Kupavskii,
	The Erd\H{o}s matching conjecture and concentration inequalities,
	\emph{J. Combin. Theory Ser. B} \textbf{157} (2022), 366--400.
	
	\bibitem{FLMW}
	P. Frankl, H. Lu, J. Ma and Y. Wu,
	\emph{Towards the Erd\H{o}s matching conjecture for $4$-uniform
		hypergraphs: stability and applications},
	arXiv:2602.19230v2, 2026.
	
	\bibitem{GGMS}
	L. Gishboliner, S. Glock, P. Michaeli and A. Sgueglia,
	Defect and transference versions of the Alon--Frankl--Lov\'asz theorem,
	\emph{Combin. Probab. Comput.} \textbf{35} (2026), 542--553.
	Version with Appendix A: arXiv:2503.05089v2.
	
	\bibitem{GKM}
	L. Gishboliner, M. Krivelevich and P. Michaeli,
	Color-biased Hamilton cycles in random graphs,
	\emph{Random Structures Algorithms} \textbf{60} (2022), 289--307.
	
	\bibitem{KNRS}
	Y. Kohayakawa, B. Nagle, V. R\"odl and M. Schacht,
	Weak hypergraph regularity and linear hypergraphs,
	\emph{J. Combin. Theory Ser. B} \textbf{100} (2010), 151--160.
	
	\bibitem{LTWY}
	H. Liu, M. Turevskii, L. Wang and Z. Yan,
	\emph{Structural reductions for monochromatic matchings and Ramsey tilings},
	arXiv:2606.24863, 2026.
	
	\bibitem{Lovasz78}
	L. Lov\'asz,
	Kneser's conjecture, chromatic number, and homotopy,
	\emph{J. Combin. Theory Ser. A} \textbf{25} (1978), 319--324.
	
	\bibitem{LP}
	L. Lov\'asz and M. D. Plummer,
	\emph{Matching Theory},
	Annals of Discrete Mathematics, vol. 29,
	North-Holland, Amsterdam, 1986.
	
	\bibitem{Matousek}
	J. Matou\v{s}ek,
	\emph{Using the Borsuk--Ulam Theorem: Lectures on Topological
		Methods in Combinatorics and Geometry},
	Universitext, Springer, Berlin, 2003.
	
	\bibitem{Rockafellar}
	R. T. Rockafellar,
	\emph{Convex Analysis},
	Princeton Mathematical Series, vol. 28,
	Princeton University Press, Princeton, NJ, 1970.
	
	\bibitem{Schrijver}
	A. Schrijver,
	\emph{Combinatorial Optimization: Polyhedra and Efficiency},
	Algorithms and Combinatorics, vol. 24,
	Springer, Berlin, 2003.
	
	\bibitem{XYZ}
	C. Xu, H. Yang and S. Zhang,
	On characterizing the critical graphs for matching Ramsey numbers,
	\emph{Discrete Appl. Math.} \textbf{287} (2020), 15--20.
	
\end{thebibliography}
\end{document}